\documentclass[english,11pt]{amsart}
\usepackage{amsopn,amsthm,amsfonts,amsmath,amssymb,amscd}
\usepackage{babel}
\usepackage{amssymb,tikz}

\newtheorem{Theorem}{Theorem}[section]
\newtheorem{Lemma}[Theorem]{Lemma}

\newtheorem{Corollary}[Theorem]{Corollary}

\newtheorem*{Remark}{Remark}
\newtheorem{Def}[Theorem]{Definition}

\newenvironment{Proof}{{\bf Proof.}}{\hfill $\blacksquare$}
\newenvironment{Proof*}{{\it Proof.}}

\newcommand{\ZZ}{\mathbb{Z}}

\newcommand{\prob}{{\rm Prob}}

\newcommand{\rk}{\mathrm{rk}}

\newcommand{\annr}{\mathrm{Ann}_r}

\begin{document}

\title{The multiplication probability of a finite ring}

\author{David Dol\v zan}

\address{D.~Dol\v zan:~Department of Mathematics, Faculty of Mathematics
and Physics, University of Ljubljana, Jadranska 19, SI-1000 Ljubljana, Slovenia, and Institute of Mathematics, Physics and Mechanics, Jadranska 19, SI-1000 Ljubljana, Slovenia; e-mail: 
david.dolzan@fmf.uni-lj.si}

\subjclass[2020]{16P10, 16L30, 16K20} 
\keywords{finite ring; Jacobson radical; probability}
\thanks{The author acknowledges the financial support from the Slovenian Research Agency  (research core funding No. P1-0222)}

\begin{abstract} 
We study the probability that the product of two randomly chosen elements in a finite ring $R$ is equal to some fixed element $x \in R$. We calculate this probability for semisimple rings and some special classes of local rings, and find the bounds for this probability for an arbitrary finite ring.
\end{abstract}

\maketitle 

 \section{Introduction}

\bigskip

Problems concerning discrete probabilities have been historically first studied in the setting of groups, see for example \cite{barry, guralnick, gustafson, lescot, rusin}. 
Later, studies in this vein have also been carried out by researchers studying rings. Authors studied the commuting probability (also called commutativity degree)  (see \cite{basnet, buckley1, buckley, dutta, dutta1, machale, shumyatsky}) and the probability of zero multiplication (also called the nullity degree) (see \cite{dolzan, esmkhani, esmkhani1, salih, shumyatsky}) of finite rings. Recently, the authors in \cite{rehman} and \cite{rehman0} tackled the problem of finding the multiplication probability of a finite ring, defined by $\prob_x(R)=\frac{|\{(a,b) \in R^2; ab=x\}|}{|R|^2}$ where $x \in R$ is an arbitrary element. Obviously, this is a generalization of the nullity degree, which is why this probability is sometimes also called the \emph{generalized probability} of a finite ring.
However, the authors in \cite{rehman, rehman0} have tackled this problem only in the setting of commutative  rings and they mainly concerned themselves with the ring of integers modulo $n$ for some fixed integer $n$.

In this paper, we investigate the generalized probability of an arbitrary finite ring. 
In the following section, we gather some definitions and preliminary results that will come in handy throughout the paper. We also establish the bounds for the multiplication probability (see Lemma \ref{unit}), thereby generalizing the statements of  \cite[Theorem 4.1]{dolzan} and  \cite[Theorem 2.1]{rehman}. We also calculate the multiplication probability of a direct product of rings and find the connection between the multiplication probabilities of a ring and its factor ring over the Jacobson radical. In Section 3, we find the multiplication probability of a simple ring (see Theorem \ref{matrix}). In the final section, we study the multiplication probability of local rings. Firstly, we establish the bounds and examine when these bounds are achieved (see Theorem \ref{unitlocal} and  Corollary \ref{equivalence}). Then we calculate the multiplication probability in two special classes of finite local rings: rings such that the Jacobson radical $J$ has the maximal possible nilpotency index (see Theorem \ref{veriga}) - such classes of rings contain for example all rings $\ZZ_n$ where $n$ is a power of a prime, so Theorem \ref{veriga} is a generalization of \cite[Theorem 2.1 and Theorem 2.2]{rehman0}, and also the rings where $J^2=0$ (see Theorem \ref{kvadratnula}).

\bigskip

 \section{Definitions and preliminaries}

Throughout the paper, $R$ will denote a finite ring with unity $1$. We will denote its group of units by $R^*$, its set of zero-divisors by $Z(R)=\{x \in R; \text{ there exists } 0 \neq y \in R \text{ such that } xy=0 \text{ or } yx=0\}$ and its Jacobson radical by $J(R)$. When the ring in question will be beyond doubt, we will write simply $J$ instead of $J(R)$. We shall often denote the equivalence class of $x$ in $R/J$ by $\overline{x}=x+J$. The (right) annihilator of an element $x \in R$ will be denoted by $\annr(x)=\{y \in R; xy=0\}$. If the ring in question will not be obvious, we shall write $\annr(x, R)$ instead of $\annr(x)$.


We will denote by $M_n(F)$ the ring of $n$-by-$n$ matrices over the field $F$ and we will denote the rank of $A \in M_n(F)$ by $\rk(A)$.



We shall make use of the following simple lemma, which shows that in a finite ring, every left zero-divisor is also a right zero-divisor (and vice-versa).

\begin{Lemma}\cite[Proposition 1.3]{kobayashi}
\label{leftright}
    Let $R$ be a finite ring and $a \in R$. If there exists $0 \neq b \in R$ such that $ba=0$ then there exists $0 \neq c \in R$ such that $ac=0$. 
\end{Lemma}

\bigskip

In order to calculate the generalized probability, we have the following definition.

\begin{Def}
    Let $a, x \in R$. We define 
    $$\delta_x(a)=\begin{cases} 1, \text { if there exists } b \in R \text { such that } ab=x, \\
    0, \text { otherwise}.
    \end{cases}$$
\end{Def}

\bigskip

We immediately have the following lemma.

\begin{Lemma}
\label{sum}
Let $R$ be a finite ring. 
    For $x \in R$, we have $\prob_x(R)=\frac{\sum_{a \in R}|\annr(a)|\delta_x(a)}{|R|^2}$.
\end{Lemma}
\begin{Proof}
    Assume that there exists $b \in R$ such that $ab=x$. Now, $ac=x$ implies $a(c-b)=0$, so $c-b \in \annr(a)$. Thus, $ac=x$ for exactly all $c \in b + \annr(a)$.
\end{Proof}

\bigskip

The following lemma is a generalization of \cite[Theorem 4.1]{dolzan} and \cite[Theorem 2.1]{rehman}. Note that the fact that $\prob_0(R) \leq \frac{2|R|-2|Z(R)|+|Z(R)|^2}{|R|^2}$ was already proved in \cite[Theorem 4.1]{dolzan}.

\begin{Lemma}
\label{unit}
    If $R$ is a finite ring then $x \in R^*$ if and only if $\prob_x(R)=\frac{|R^*|}{|R|^2}$.
    Furthermore, if $0 \neq x \notin R^*$ then $\frac{|R^*|}{|R|^2} \cdot \frac{2+|Z(R)|}{|Z(R)|} \leq \prob_x(R) \leq \frac{|R|-2|Z(R)|+|Z(R)|^2}{|R|^2}$ and $\frac{2|R|+|Z(R)|-2}{|R|^2} \leq \prob_0(R) \leq \frac{2|R|-2|Z(R)|+|Z(R)|^2}{|R|^2}$.
\end{Lemma}
\begin{Proof}
    Denote $A=\{(a,b) \in R^2; ab=x\}$.
    
    Suppose firstly that $x \in R^*$.
    Since, $ab = x$ if and only if $x^{-1}ab=1$, we see that $(a,b) \in A$ if and only if $b \in R^*$ and $(xb^{-1},b) \in A$, so $|A|=|R^*|$ and thus $\prob_x(R)=\frac{|A|}{|R|^2}=\frac{|R^*|}{|R|^2}$.
    
    If $x=0$ then $(0,b), (a,0) \in A$ for every $b \in R$ and every $a \in R$. However, if $0 \neq a \in Z(R)$, then $|\annr(a)| \geq 2$ by Lemma \ref{leftright}, so $|A| \geq |R| + |R^*| + 2(|Z(R)|-1)=2|R|+|Z(R)|-2$ and the lower bound follows. The upper bound follows from \cite[Theorem 4.1]{dolzan}. Note also that in this case, $|A| > |R^*|$ so $\prob_x(R) \neq \frac{|R^*|}{|R|^2}$.
    
    Finally, suppose $0 \neq x \in Z(R)$. If $a \in R^*$ then $(a,a^{-1}x), (xa,a^{-1}+t) \in A$ for every $t \in \annr(xa)$. Since $xa \in Z(R)$, we have by Lemma \ref{leftright} that $|\annr(xa)| \geq 2$.
    If $xa=xb$ for some $a,b \in R^*$ then $x(a-b)=0$, so $a-b \in Z(R)$. Thus there exists $z \in Z(R)$ such that $a=b+z$. This implies that we have at least $\frac{|R^*|}{|Z(R)|}$ different elements $xa$, so $|A| \geq |R^*| + \frac{2|R^*|}{|Z(R)|}$ and thus $\prob_x(R) \geq \frac{|R^*|}{|R|^2} \cdot \frac{2+|Z(R)|}{|Z(R)|}$. Observe that again in this case, $|A| > |R^*|$ so $\prob_x(R) \neq \frac{|R^*|}{|R|^2}$. Since every element of $R$ is either a unit or a zero-divisor, this now proves that $x \in R^*$ if and only if $\prob_x(R)=\frac{|R^*|}{|R|^2}$.
    To prove the upper bound, observe that
    for $a \in R^*$, we have $|\annr(a)|=1$ and $\delta_x(a)=1$, $\delta_x(0)=0$ and for $0 \neq a \in Z(R)$, we have $\delta_x(a) \leq 1$ and $|\annr(a)| \leq |Z(R)|$, so $|A| \leq |R^*| + |Z(R)|(|Z(R)|-1)$.
\end{Proof}

\bigskip

The following lemma shows that we may always limit ourselves to studying the multiplication probability of directly indecomposable rings.

\begin{Lemma}
    \label{product}
    If $R=R_1 \times R_2 \times \ldots \times R_n$ is a direct product of finite rings and $x=(x_1,x_2, \ldots, x_n) \in R$ then $\prob_x(R)=\prod_{i=1}^n\prob_{x_i}(R_i)$.
\end{Lemma}
\begin{Proof}
    Observe that for any elements $a=(a_1,a_2,\ldots,a_n), b=(b_1,b_2,\ldots,b_n)$ and $x=(x_1,x_2,\ldots,x_n) \in R$ we have $ab=x$ if and only if $a_ib_i=x_i$ for every $i=1,2,\ldots,n$.
\end{Proof}

\bigskip

We conclude this preliminary section with the following estimate.

\begin{Lemma}
    \label{overj}
Let $R$ be a finite ring and $I$ a proper ideal in $R$. 
    For any $x \in R$, we have $\prob_x(R) \leq \prob_{{x+I}}(R/I)$.
\end{Lemma}
\begin{Proof}
  We have $\prob_x(R)=\frac{\sum_{a \in R}|\annr(a)|\delta_x(a)}{|R|^2}$ by Lemma \ref{sum}.
  Since $\delta_x(a)=1$ implies $\delta_{{x+I}}({a+I})=1$ and since for any $|I|$ elements in $\annr(a,R)$ we have at least one element in $\annr({a+I},R/I)$, we get
  $\prob_x(R) \leq \frac{\sum_{a \in R}|\annr({a+I},R/I)|I||\delta_{{x+I}}(a+I)}{|R|^2}=\frac{\sum_{{a+I}\in R/I}|\annr({a+I},R/I)|I|^2|\delta_{{x+I}}({a+I})}{|R|^2}=\prob_{{x+I}}(R/I)$.
\end{Proof}

\bigskip

 \section{Simple rings}

\bigskip

In this section, we investigate the multiplication probability of an arbitrary (semi)simple finite ring. Lemma \ref{product} implies that to calculate the multiplication probability of a semisimple ring, we have to study the multiplication probability of a simple ring, so we limit ourselves to this special case.

Since the $k$-dimensional subspaces in $F^n$ that contain a certain subspace $U$ are in one to one corresponence with $(k-r)$-dimensional subspaces in the quotient space $F^n/U \simeq F^{n-r}$, we immediately have the following lemma.

\begin{Lemma}\cite[Theorem 13.1]{andrews}
\label{subspaces}
Let $F$ be a finite field with $q$ elements and let $r \leq k \leq n$ be integers. If $U$ is a $r$-dimensional subspace in $F^n$, then there exist exactly $\prod_{i=0}^{k-r-1}\frac{q^{n-r}-q^i}{q^{k-r}-q^i}$ $k$-dimensional subspaces in $F^n$ that contain $U$. \end{Lemma}
  

\bigskip

We now have the following theorem.

\begin{Theorem}
\label{matrix}
Let $F$ be a finite field with $q$ elements and $n$ an integer. For a matrix $X \in M_n(F)$ of rank $r \leq n$, we have $$\prob_X(M_n(F))=\sum_{k=r}^n\frac{ \prod_{i=0}^{k-r-1}\frac{q^{n-r}-q^i}{q^{k-r}-q^i} \prod_{i=0}^{k-1}(q^n-q^i) }{q^{n(n+k)}}.$$
\end{Theorem}
\begin{Proof}
    By Lemma \ref{sum}, we have $\prob_X(R)=\frac{\sum_{k=0}^n\sum_{\rk(A)=k}|\annr(A)|\delta_X(A)}{|R|^2}$.
    Note that there exists $B \in M_n(F)$ such that $AB=X$ if and only if the column space of $X$ is a subspace of the column space of $A$. Denote the column space of $X$ by $U$. Since $\rk(X)=r$, we have to count all matrices $A$ of rank $k \geq r$ such that their column space contains $U$. By Lemma \ref{subspaces} there exist exactly $\prod_{i=0}^{k-r-1}\frac{q^{n-r}-q^i}{q^{k-r}-q^i}$ $k$-dimensional subspaces in $F^n$ that contain $U$.
    Choose any such subspace $V$. Now, we have as many matrices with their column space equal to $V$ as there are surjective mappings from $F^n$ to $V$, which is exactly $(q^n-1)(q^n-q)\ldots(q^n-q^{k-1})$. Thus, there are exactly $\prod_{i=0}^{k-r-1}\frac{q^{n-r}-q^i}{q^{k-r}-q^i} \prod_{i=0}^{k-1}(q^n-q^i)$ matrices $A$ of rank $k$ such that $\delta_X(A)=1$. Choose $x \in F^n$. Since any matrix $A$ of rank $k$ such that the equation $Ay=x$ is solvable has exactly $q^{n-k}$ solutions $y \in F^n$, we see that there exist exactly $(q^{n-k})^n$ matrices $B$ such that $AB=X$ for any $A \in M_n(F)$ such that $\delta_X(A)=1$.
    This implies that $$\prob_X(M_n(F))=\frac{\sum_{k=r}^n q^{n(n-k)}\prod_{i=0}^{k-r-1}\frac{q^{n-r}-q^i}{q^{k-r}-q^i} \prod_{i=0}^{k-1}(q^n-q^i) }{q^{2n^2}},$$ which proves the theorem.
\end{Proof}

\bigskip

 \section{Local rings}

\bigskip

In this section, we shall investigate the multiplication probability of local finite rings.
We will need the following lemma.

\begin{Lemma}\cite[Theorem 2]{raghavendran}
\label{vsejeq}
    Let $R$ be a finite local ring with $|R|=q^n$ and $R/J=GF(q)$. Let $I$ be a left (right) ideal in $R$. Then there exists integer $\alpha$ such that $|I|=q^\alpha$.
    Furthermore, there exists $g \in R^*$ such that $g+J$ is a primitive element of the field $R/J$ and $g^k-g^l \in J$ for some $k,l \in \{0,1,\ldots,q-1\}$ implies $k=l$.
\end{Lemma}

\bigskip

From Lemma \ref{unit}, we now get the following.

\begin{Theorem}
\label{unitlocal}
    Let $R$ be a finite local ring with $|R|=q^n$ for some $n \geq 2$ and $R/J=GF(q)$.
    Then $x \in R^*$ if and only if $\prob_x(R)=\frac{q-1}{q^{n+1}}$.
    Furthermore, if $0 \neq x \in J$ then $\frac{(q-1)(q^{ n-2}+1)}{q^{2n-1}} \leq \prob_x(R) \leq \frac{q^{ n-1}+q-2}{q^{n+1}}$ and $\frac{3q^{n-1}-q^{n-2}-1}{q^{2n-1}} \leq \prob_0(R) \leq \frac{q^{ n-1}+2q-2}{q^{n+1}}$.
\end{Theorem}
\begin{Proof}
    For $x \in R^*$ the statement follows directly from Lemma \ref{unit}. 
    
    If $x=0$ then $(0,b), (a,0) \in A$ for every $b \in R$ and every $a \in R^*$. However, if $0 \neq a \in J$, then $\annr(a) \neq \{0\}$ by Lemma \ref{leftright}, so $|\annr(a)| \geq q$ by Lemma \ref{vsejeq}. Thus $|A| \geq |R| + |R^*| + q(|J|-1)=2|R|+(q-1)|J|-q$, so $\prob_0(R) \geq \frac{3q^{n-1}-q^{n-2}-1}{q^{2n-1}}$. The upper bound follows from Lemma \ref{unit} and the fact that $|Z(R)|=|J|=q^{n-1}$.
    
    Finally, suppose $0 \neq x \in J$. If $a \in R^*$ then $(a,a^{-1}x), (xa,a^{-1}+t) \in A$ for every $t \in \annr(xa)$. Since $\annr(xa)$ is nontrivial, we have $|\annr(xa)| \geq q$ by Lemma \ref{vsejeq}. Now, if $xa=xb$ for some $a,b \in R^*$ then $x(a-b)=0$, so $a-b \in J$. Thus there exist at least $q-1$ different elements $xa$, so $|A| \geq |R^*| + q(q-1)$ and thus $\prob_x(R) \geq \frac{(q-1)(q^{ n-1}+q)}{q^{2n}}=\frac{(q-1)(q^{ n-2}+1)}{q^{2n-1}}$.  The upper bound again follows directly from Lemma \ref{unit}.
\end{Proof}

\bigskip

The next two corollaries show us when exactly these bounds are achieved.

\begin{Corollary}
\label{equivalence}
    Let $R$ be a finite local ring with $|R|=q^n$ for some $n \geq 2$ and $R/J=GF(q)$. The following statements are equivalent.
    \begin{enumerate}
        \item $\prob_x(R) = \frac{(q-1)(q^{ n-2}+1)}{q^{2n-1}}$ for every $0 \neq x \in J$.
        \item $\prob_x(R) = \frac{q^{ n-1}+q-2}{q^{n+1}}$ for every $0 \neq x \in J$.
        
        \item
        $\prob_0(R) = \frac{3q^{n-1}-q^{n-2}-1}{q^{2n-1}}$.


        \item 
        $R$ is a ring with $q^2$ elements.
    \end{enumerate}
\end{Corollary}
\begin{Proof}
   (1) $\Rightarrow$ (4): By the proof of Theorem \ref{unitlocal}, we see that $|\annr(x)|=q$ and that $x=ab$ implies that either $a \in R^*$ or $b \in R^*$. Since this is true for any $x \in J$, we have $J^2=0$. The fact that $|\annr(x)|=q$ now implies $|J|=q$, so $|R|=q^2$.

   (2) $\Rightarrow$ (4): The proof of Lemma \ref{unit} shows that for any $0 \neq a \in J$, we have $|\annr(a)| = |J|$, which implies that $J^2=0$. We also have $\delta_x(a)=1$ for every $0 \neq a \in J$, which implies that for every $0 \neq a \in J$ there exists $u_a \in R^*$ such that $au_a=x$, so $a=xu_a^{-1}$. But if $u_1 - u_2 \in J$ then $xu_1=xu_2$, so $|xR| \leq q$, which implies $|J|=q$, and therefore $|R|=q^2$.
   
   (3) $\Rightarrow$ (4): We know that $|\annr(a)| = q$ for every $0 \neq a \in J$ by the proof of Theorem \ref{unitlocal}. Suppose $J^k \neq 0$ and $J^{k+1}=0$ for some integer $k$. Then $|J^k| \geq q$ by Lemma \ref{vsejeq}. If $k \geq 2$, then $|J| > q^2$, so there exists $a \in J^k$ such that $|\annr(a)| \geq q^2$, a contradiction. We have therefore proved that $J^2=0$. Since $|\annr(a)| = q$ for every $a \in J$ and $aJ=0$ this  implies that $|J|=q$, therefore $|R|=q^2$.

   (4) $\Rightarrow$ (1), (2), (3):
   If $|R|=q^2$, then $|J|=q$ and $J^2=0$, so $|\annr(a)| = q$ for every $0 \neq a \in J$. Now, it is easy now check that all inequalities in the proof of Corollary \ref{unitlocal} are actually equalities, so (1) and (3) follow directly. 
   Observe also that the statements (1) and (2) are equal when $n=2$, so (2) follows as well.
\end{Proof}

\begin{Corollary}
\label{finalbound}
   Let $R$ be a finite local ring with $|R|=q^n$ and $R/J=GF(q)$. Then $\prob_0(R)=       \frac{q^{n-1}+2q-2}{q^{n+1}}$ if and only if $J^2=0$.
\end{Corollary}
\begin{Proof}
  This follows directly from \cite[Theorem 4.1]{dolzan}.
\end{Proof}

\bigskip

In some special cases of local rings, we can actually calculate all the multiplication probabilities. We shall need the following lemma.

\begin{Lemma}
\label{inver}
    Let $R$ be a finite ring and let $u \in R^*$ and $j \in J$. Then $u+j \in R^*$.
\end{Lemma}
\begin{Proof}
  Note that $u+j=u(1+u^{-1}j)$ and $j'=u^{-1}j \in J$. Since every element of $J$ is nilpotent, we have $j'^n=0$ for some integer $n$, therefore $(1+j')^{-1}=1-j'+(j')^2+\ldots+(-1)^{n-1}(j')^{n-1}$.
\end{Proof}

\bigskip

We now have the following theorem, which generalizes Theorems 2.1 and 2.2 from \cite{rehman0}.

\begin{Theorem}
\label{veriga}
    Let $R$ be a finite local ring with $|R|=q^n$ and $R/J=GF(q)$
    and suppose that $J^{n-1} \neq 0$. If $x \in J^k \setminus J^{k+1}$ for some integer $k$, then $\prob_x(R)=\frac{(k+1)(q-1)}{q^{n+1}}$. Furthermore, $\prob_0(R)=\frac{(n+1)q-n}{q^{n+1}}$.
\end{Theorem}
\begin{Proof}
   If $J^i=J^{i+1}$ for some $i$, then $J^i=0$ by the Nakayama's Lemma. Therefore $J \supset J^2 \supset \ldots \supset J^{n-1} \supset 0$ is a strictly decreasing chain of ideals. Denote $J^0=R$ and observe that for every $i=0,1,\ldots,n-1$, we have that $J^i/J^{i+1}$ is a (nontrivial) left $R/J$ vector space with the scalar multiplication defined by $(x+J)(j+J^{i+1})=xj+J^{i+1}$. This implies that $|J^i|=|J^{i+1}|q^{k_i}$ for some positive integer $k_i$ and therefore $q^n=|R|=q^{k_0+k_2+\ldots+k_{n-1}}$, which of course implies that $k_0=k_1=\ldots=k_{n-1}=1$. This shows that $|J^i|=q^{n-i}$ for every $i=0,1,\ldots,n$.
   Let $g \in R^*$ be the the element from Lemma \ref{vsejeq}.
   Now, choose $x_i \in J^i \setminus J^{i+1}$ for every $i \in \{1,2,\ldots,n-1\}$. 
   Choose also any $k \in \{0,1,\ldots,n-1\}$ and observe that $X_k=\{\sum_{i=k}^{n-1}\lambda_i x_i;  \lambda_i \in \{0, g, g^2, \ldots, g^{q-1} \} \text { for every } i \in \{k,k+1,\ldots,n-1\}\} \subseteq J^k$. Let us prove by induction that \begin{equation}
   \label{eqjnak}
       \sum_{i=k}^{n-1}\lambda_ix_i=\sum_{i=k}^{n-1}\mu_ix_i
   \end{equation} where $\lambda_i, \mu_i \in \{0, g, g^2, \ldots, g^{q-1} \}$ for every $i \in \{k,k+1,\ldots,n-1\}$ implies $\lambda_i=\mu_i$ for every $i \in \{k, k+1,\ldots,n-1\}$. Since $(\lambda_k - \mu_k)x_k \in J^{k+1}$ and $x_k \notin J^{k+1}$, we have either $\lambda_k=\mu_k=0$ or $\lambda_k=g^\alpha_k$ and $\mu_k =g^\beta_k$ for some $\alpha_k, \beta_k \in \{0,1,\ldots,q-1\}$. By the properties of the element $g$, this implies $\alpha_k=\beta_k$, so $\lambda_k=\mu_k$. Now, assume that $\lambda_i = \mu_i$ for all $i=k,k+1,\ldots,j$. Then Equation (\ref{eqjnak}) gives us $\sum_{i=j+1}^{n-1}\lambda_i x_i=\sum_{i=j+1}^{n-1}\mu_i x_i$ and we can use the same argument as above to establish that $\lambda_{j+1}=\mu_{j+1}$. 
   
   Since $|J^k|=q^{n-k}$, we see that $J^k=X_k$. Similar argument also gives us $J^k=Y_k=\{\sum_{i=k}^{n-1}x_i \mu_i;  \mu_i \in \{0, g, g^2, \ldots, g^{q-1} \} \text { for every } i \in \{k,k+1,\ldots,n-1\}\}$. Since $x_1g \in J \setminus J^2$, we see that $x_1g=g^\alpha x_1 + j$ for some $\alpha \in \{0,1,\ldots,q-1\}$ and some $j \in J^2$. The facts that $J=X_1$, $J^n=0$ and $J^{n-1} \neq 0$ now imply that $x_1^{n-1} \neq 0$, so $x_1^i \in J^i \setminus J^{i+1}$ for every $i=1,2,\ldots,n-1$. 
   
   Now, choose any $x \in J^k \setminus J^{k+1}$. By the above, we have $x=\sum_{i=k}^{n-1}\lambda_i x_1^i$ for some $\lambda_i \in \{0, g, g^2, \ldots, g^{q-1} \}$ with $\lambda_k \neq 0$. Choose also any $y=\sum_{i=k}^{n-1}x_1^i \mu_i \in J^k$, with $\mu_i \in \{0, g, g^2, \ldots, g^{q-1} \}$. Note that $x'=\sum_{i=k}^{n-1}\lambda_i x_1^{i-k} \in R^*$ by Lemma \ref{inver}, so $xx'^{-1}=x_1^k$ and therefore $xx'^{-1}(\sum_{i=k}^{n-1}x_1^{i-k} \mu_i)=y$.   
   We have therefore proved that $xR = J^k$.  Since $xR \simeq R/\annr(x)$ this means that $|\annr(x)|=q^{k}$ and because $J^{n-k} \subseteq \annr(x)$, we have proved that $\annr(x)=J^{n-k}$. Choose $a \in R$ and note that $\delta_x(a)=0$ if and only if $a \in J^{k+1}$.
   Now, Lemma \ref{sum} gives us $\prob_x(R)=\frac{\sum_{a \in R \setminus J^{k+1}}|\annr(a)|}{|R|^2}=\frac{\sum_{i=0}^{k}\sum_{a \in J^i \setminus J^{i+1}}|\annr(a)|}{|R|^2}$. For any $a \in J^i \setminus J^{i+1}$, we have $\annr(a)=J^{n-i}$, so $|\annr(a)|=q^i$. Furthermore, $|J^i \setminus J^{i+1}|=q^{n-i}-q^{n-i-1}$, which gives us $\prob_x(R)=\frac{\sum_{i=0}^{k}{q^{n-1}(q-1)}}{|R|^2}=\frac{(k+1)(q-1)}{q^{n+1}}$.

   Finally, if $x=0$, then $\delta_x(a)=1$ for every $a \in R$, so we have that $\prob_x(R)=\frac{|R|+ \sum_{i=0}^{n-1}\sum_{a \in J^i \setminus J^{i+1}}|\annr(a)|}{|R|^2}=\frac{q^n+\sum_{i=0}^{n}q^{n-1}(q-1)}{q^{2n}}=\frac{q+n(q-1)}{q^{n+1}}=\frac{(n+1)q-n}{q^{n+1}}$.
   \end{Proof}

\begin{Remark}
   Observe that Theorem \ref{veriga} gives us the direct explicit value of the $\prob_x(R)$ for any $x \in R$, in the case $R=\ZZ_{p^r}$ for some integer $r$ and prime $p$. However, Theorems 2.1 and 2.2 from \cite{rehman0} only give us the summation formula. Also, Lemma \ref{veriga} together with Lemma \ref{product} now also give us a formula for $\prob_x(R)$ in case $R=\ZZ_{n}$ for any integer $n$ and any $x \in R$, which generalizes Theorems 2.10, 2.11, 2.12, 2.13, 2.14 and 2.15 from \cite{rehman}.
\end{Remark}

\bigskip

We also have the following theorem.

\begin{Theorem}
\label{kvadratnula}
    Let $R$ be a finite local ring with $|R|=q^n$ and $R/J=GF(q)$
    and suppose that $J^2 = 0$. Then $\prob_x(R)=\begin{cases}
        \frac{2(q-1)}{q^{n+1}}, \text {if } 0 \neq x \in J, \\
        \frac{q-1}{q^{n+1}}, \text {if } x \notin J, \\
        \frac{q^{n-1}+2q-2}{q^{n+1}}, \text { if } x=0.
    \end{cases}$
\end{Theorem}
\begin{Proof}
  Choose $0 \neq x \in J$ and then choose $0 \neq a \in J$ such that $\delta_x(a)=1$. Since $J^2=0$, if there exists $b \in R$ such that $ab=x$, we know that $b \notin J$, so $b$ is invertible and thus $a=xb^{-1}$. However, if $xb_1=xb_2$ for two different invertible elements $b_1$ and $b_2$, then $x(b_1-b_2)=0$, so $b_1-b_2 \in J$ and therefore there exist at least $q-1$ different non-zero elements of the form $xb$ for some invertible $b \in R$. Since $J^2=0$, we have $|xR \setminus \{0\}|=q-1$, so we have exactly $q-1$ different elements $a \in J$ such that $\delta_x(a)=1$. For each such $a$, we have $\annr(a)=J$. On the other hand, if $a \notin J$, then $a$ is invertible, so $\delta_x(a)=1$ and $\annr(a)=\{0\}$.
  Therefore by Lemma \ref{sum}, we have $\prob_x(R)=\frac{(q-1)|J|+|R^*|}{|R|^2}=\frac{2(q-1)q^{n-1}}{q^{2n}}=\frac{2(q-1)}{q^{n+1}}$.
  
  If $x \notin J$, we have $\prob_x(R)=\frac{q-1}{q^{n+1}}$ by Theorem \ref{unitlocal}. 
  
  Finally, if $x=0$, we have $|\annr(a)|=1$ for invertible elements $a$, $|\annr(a)|=q^{n-1}$ for non-invertible non-zero elements $a$ and of course $|\annr(0)|=q^n$, so $\prob_x(R)=\frac{q^{n-1}(q-1)+(q^{n-1}-1)q^{n-1}+q^n}{q^{2n}}=\frac{q^{n-1}+2q-2}{q^{n+1}}$.
\end{Proof}

\bigskip

{\bf Acknowledgement} \\

The author would like to thank the anonymous referee who provided useful and detailed comments that helped to improve the quality of this paper.

\bigskip

{\bf Statements and Declarations} \\

The author states that there are no competing interests. 

\bigskip

\bibliographystyle{amsplain}
\bibliography{biblio}

\providecommand{\bysame}{\leavevmode\hbox to3em{\hrulefill}\thinspace}
\providecommand{\MR}{\relax\ifhmode\unskip\space\fi MR }
\providecommand{\MRhref}[2]{%
  \href{http://www.ams.org/mathscinet-getitem?mr=#1}{#2}
}
\providecommand{\href}[2]{#2}
\begin{thebibliography}{10}

\bibitem{andrews}
George~E. Andrews, \emph{The theory of partitions}, Cambridge Mathematical Library, Cambridge University Press, Cambridge, 1998, Reprint of the 1976 original. \MR{1634067}

\bibitem{barry}
F.~Barry, D.~MacHale, and \'{A}. N\'{\i}~Sh\'{e}, \emph{Some supersolvability conditions for finite groups}, Math. Proc. R. Ir. Acad. \textbf{106A} (2006), no.~2, 163--177. \MR{2266824}

\bibitem{basnet}
Dhiren~Kumar Basnet and Jutirekha Dutta, \emph{Some bounds for commuting probability of finite rings}, Proc. Indian Acad. Sci. Math. Sci. \textbf{129} (2019), no.~1, Paper No. 1, 6. \MR{3887207}

\bibitem{buckley1}
S.~M. Buckley, D.~MacHale, and Y.~Zelenyuk, \emph{Finite rings with large anticommuting probability}, Appl. Math. Inf. Sci. \textbf{8} (2014), no.~1, 13--25. \MR{3117772}

\bibitem{buckley}
Stephen~M. Buckley and Desmond MacHale, \emph{Commuting probability for subrings and quotient rings}, J. Algebra Comb. Discrete Struct. Appl. \textbf{4} (2017), no.~2, 189--196. \MR{3601350}

\bibitem{dolzan}
David Dol\v{z}an, \emph{The probability of zero multiplication in finite rings}, Bull. Aust. Math. Soc. \textbf{106} (2022), no.~1, 83--88. \MR{4448946}

\bibitem{dutta}
Jutirekha Dutta, Dhiren~Kumar Basnet, and Rajat~Kanti Nath, \emph{On commuting probability of finite rings}, Indag. Math. (N.S.) \textbf{28} (2017), no.~2, 372--382. \MR{3624561}

\bibitem{dutta1}
\bysame, \emph{Characterization and commuting probability of {$n$}-centralizer finite rings}, Proyecciones \textbf{42} (2023), no.~6, 1489--1498. \MR{4679483}

\bibitem{esmkhani}
M.~A. Esmkhani and S.~M. Jafarian~Amiri, \emph{The probability that the multiplication of two ring elements is zero}, J. Algebra Appl. \textbf{17} (2018), no.~3, 1850054, 9. \MR{3760023}

\bibitem{esmkhani1}
\bysame, \emph{Characterization of rings with nullity degree at least {$\frac14$}}, J. Algebra Appl. \textbf{18} (2019), no.~4, 1950076, 8. \MR{3928520}

\bibitem{guralnick}
Robert~M. Guralnick and Geoffrey~R. Robinson, \emph{On the commuting probability in finite groups}, J. Algebra \textbf{300} (2006), no.~2, 509--528. \MR{2228209}

\bibitem{gustafson}
W.~H. Gustafson, \emph{What is the probability that two group elements commute?}, Amer. Math. Monthly \textbf{80} (1973), 1031--1034. \MR{327901}

\bibitem{kobayashi}
Yuji Kobayashi and Kwangil Koh, \emph{A classification of finite rings by zero divisors}, J. Pure Appl. Algebra \textbf{40} (1986), no.~2, 135--147. \MR{830317}

\bibitem{lescot}
Paul Lescot, \emph{Isoclinism classes and commutativity degrees of finite groups}, J. Algebra \textbf{177} (1995), no.~3, 847--869. \MR{1358488}

\bibitem{machale}
Desmond MacHale, \emph{Commutativity in finite rings}, Amer. Math. Monthly \textbf{83} (1975), no.~1, 30--32. \MR{384861}

\bibitem{salih}
Haval~M. Mohammed~Salih, \emph{On the probability of zero divisor elements in group rings}, Int. J. Group Theory \textbf{11} (2022), no.~4, 253--257. \MR{4509600}

\bibitem{raghavendran}
R.~Raghavendran, \emph{Finite associative rings}, Compositio Math. \textbf{21} (1969), 195--229. \MR{246905}

\bibitem{rehman}
Shafiq~ur Rehman and Muhammad~Naveed Shaheryar, \emph{On generalized probability in finite commutative rings}, Int. Electron. J. Algebra \textbf{33} (2023), 125--132. \MR{4563792}

\bibitem{rusin}
David~J. Rusin, \emph{What is the probability that two elements of a finite group commute?}, Pacific J. Math. \textbf{82} (1979), no.~1, 237--247. \MR{549847}

\bibitem{shumyatsky}
Pavel Shumyatsky and Matteo Vannacci, \emph{Commuting and product-zero probability in finite rings}, Internat. J. Algebra Comput. \textbf{34} (2024), no.~2, 201--206. \MR{4727826}

\bibitem{rehman0}
Shafiq ur~Rehman, Abdul~Qudair Baig, and Kamran Haider, \emph{A probabilistic approach toward finite commutative rings}, Southeast Asian Bull. Math. \textbf{43} (2019), no.~3, 413--418. \MR{3966299}

\end{thebibliography}

\bigskip

\end{document}